\documentclass[12pt,twoside,a4paper]{amsart}
\usepackage{amscd}
\usepackage{amsmath}
\usepackage{amssymb}
\usepackage{latexsym}   
\usepackage[english]{babel}

\input{xy}
\xyoption{all}

\newcommand{\FQSym}{\mathbf{FQSym}}

\newcommand{\QSym}{\mathbf{QSym}}
\newcommand{\WQSym}{\mathbf{WQSym}}
\newcommand{\IQSym}{\mathbf{IQSym}}
\newcommand{\SIQSym}{\mathbf{SIQSym}}

\newcommand{\Q}{\mathbf{Q}}

\newcommand{\R}{\mathcal{R}}
\newcommand{\F}{\mathcal{F}}

\newcommand{\N}{\mathbb{N}}

\newcommand{\Di}{\mathbf{Diff}}

\def\NN{{\bf N}}

\def\Sh{{\rm Sh}}
\def\QSh{{\rm QSh}}
\def\Des{{\rm Des}}
\def\desc{{\rm desc}}
\def\id{{\rm id}}
\def\Id{{\rm Id}}
\def\Hom{{\rm Hom}}
\def\lin{{\rm lin}}
\def\Lin{{\rm Lin}}
\def\coalg{{\rm coalg}}
\def\Coder{{\rm Coder}}

\def\shuff#1#2{\mathbin{
      \hbox{\vbox{\hbox{\vrule \hskip#2 \vrule height#1 width 0pt}\hrule}\vbox{\hbox{\vrule \hskip#2 \vrule height#1 width 0pt\vrule }\hrule}}}}
\def\shuffl{{\mathchoice{\shuff{5pt}{3.5pt}}{\shuff{5pt}{3.5pt}}{\shuff{3pt}{2.6pt}}{\shuff{3pt}{2.6pt}}}}
\def\shuffle{\, \shuffl \,}

\def\qshuffl{{\mathchoice{\shuff{5pt}{3.5pt}\hspace{-2.9mm}-}{\shuff{5pt}{3.5pt}\hspace{-2.9mm}-}
{\shuff{3pt}{2.6pt}\hspace{-2.2mm}-}{\shuff{3pt}{2.6pt}\hspace{-2.2mm}-}}}
\def\qshuffle{\,\qshuffl\,}

\def\<{\langle}
\def\>{\rangle}

\def\M{{\bf M}}

\def\pack{{\rm pack}}
\def\SG{{\mathfrak S}}
\def\Sym{{\bf Sym}}     % NCSF
\def\PW{{\rm PW}}

\title{Deformations of shuffles and quasi-shuffles}
\date{}
\author[L. Foissy, F. Patras, J.-Y. Thibon]
{Lo\"ic Foissy, Fr\'ed\'eric Patras,  Jean-Yves Thibon}
\address[Foissy]{Laboratoire de Math\'ematiques Pures et Appliqu\'ees Joseph Liouville\\
Universit\'e du Littoral C\^ote d’Opale, Centre Universitaire de la Mi-Voix\\
50, rue Ferdinand Buisson, CS 80699\\
62228 Calais Cedex \\ France}
\address[Patras]{Laboratoire de Math\'ematiques J.A. Dieudonn\'e\\
Universit\'e de Nice - Sophia Antipolis\\
Parc Valrose, 06108 Nice Cedex 02\\ France}
\address[Thibon]{Laboratoire d'Informatique Gaspard-Monge, Universit\'e Paris-Est, 
5, boulevard Descartes \\Champs-sur-Marne \\77454 Marne-la-Vall\'ee cedex 2 \\
France}
\email[L. Foissy]{foissy@lmpa.univ-littoral.fr}
\email[F. Patras]{patras@unice.fr}
\email[J.-Y. Thibon]{jyt@univ-mlv.fr}

\keywords{Shuffle algebras, Combinatorial Hopf algebras, Hausdorff series}

\subjclass{05E05,16T30}

\newtheorem{defi}{\indent Definition}
\newtheorem{lemma}[defi]{\indent Lemma}
\newtheorem{cor}[defi]{\indent Corollary}
\newtheorem{theo}[defi]{\indent Theorem}
\newtheorem{prop}[defi]{\indent Proposition}

\begin{document}

\begin{abstract}
We investigate  deformations of the shuffle Hopf algebra structure $\Sh(A)$
which can be defined on the tensor algebra over a commutative algebra $A$.
Such deformations, leading for example to the quasi-shuffle algebra $\QSh(A)$,
 can be interpreted as natural transformations of the functor
$\Sh$, regarded as a functor from commutative nonunital algebras  to coalgebras.
 We prove that the monoid
of natural endomophisms of the functor $\Sh$ is isomorphic to the monoid of
formal power series in one variable without constant term under composition,
so that in particular, its natural automorphisms are in bijection with
formal diffeomorphisms of the line.

These transformations can be interpreted as elements of
the Hopf algebra of word quasi-symmetric functions$\WQSym$, and in turn define deformations of its
struture. This leads to a new embedding of free quasi-symmetric functions into
$\WQSym$, whose relevance is illustrated by a simple and transparent
proof of Goldberg's formula for the coefficients of the Hausdorff series.
\end{abstract}

\maketitle

%\tableofcontents 

\section*{Introduction}
The tensor algebra over a commutative algebra $A$ is provided by the shuffle product with a commutative (Hopf) algebra structure. 
However, other products (resp. Hopf algebra) structures can be defined, 
the best known one being the quasi-shuffle product (resp. quasi-shuffle Hopf algebra).
These products arise in many contexts, {\it e.g.},  Rota--Baxter algebras, multiple zeta values (MZVs), noncommutative symmetric functions, operads ...
Moreover, they are natural: they commute with morphisms of commutative algebras; 
in other terms, they define functors from the category of commutative algebras to that of commutative Hopf algebras.

The present paper aims at studying and classifying these products and Hopf algebra structures. 
The Fa\`a di Bruno Hopf algebra plays a key role in this classification: 
its characters (in bijection with formal diffeomorphisms
tangent to the identity) happen to classify the natural deformations of shuffle algebras considered in this article.

Our approach also sheds a new light on classical constructions such as quasi-shuffles.
Most structure properties of quasi-shuffle algebras appear, from our point of view, 
as  straightforward consequences of their definition as deformations by conjugacy of shuffle algebras. 
This allows to transport automatically \it all \rm known results on shuffle algebras to quasi-shuffles 
and does not require the algebra $A$ to be graded (compare, {\it e.g.}, \cite{Hoffman}). 
We may quote the existence of a ``natural'' (but not straightforward !) gradation, and of fine
 Hopf algebraic properties such as the ones studied in \cite{npt}.

Most of these results can be phrased in terms of combinatorial properties of
Hopf algebras based on permutations 
(free quasi-symmetric functions) and surjections (word quasi-symmetric functions). 
In particular, the analysis of the relations between shuffles and quasi-shuffles leads to 
a new polynomial realization of noncommutative symmetric functions, which can be extended
to free quasi-symmetric functions. 
Beyond its naturallity from the point of view of quasi-shuffle algebras, this new realization has an interest on its own:
as an application, a simple (and, in our opinion, enlightening)
proof of Goldberg's formula for the coefficients of the Haudorff series and a generalization thereof to other, similar, series is
obtained.

\ \\
\ \\
The authors acknowledge support from the grant CARMA ANR-12-BS01-0017.

\section{The shuffle algebra over a commutative algebra}

Let $A$ be a commutative algebra over the rationals\footnote{Any field of characteristic zero would be suitable as well.}, 
(not necessarily with a unit)
and let $T(A)$ be its tensor algebra 
$$
T(A)=\bigoplus\limits_{n \in \mathbb{N}}T_n(A)=\bigoplus\limits_{n \in \mathbb{N}}A^{\otimes n}$$ 
with $A^{\otimes 0}:=\Q$.
Tensors $a_1\otimes ...\otimes a_n$ in $T_n(A)$ will be written  as words $a_1...a_n$, 
the $\otimes$ sign being reserved for tensor products of elements of $T(A)$.

The concatenation product in $T(A)$ is written $\times$ to distinguish it from the internal product in $A$
$$a_1...a_n\times b_1...b_m:=a_1...a_nb_1...b_m.$$
The (internal) product of $a$ and $b$ in $A$ will be always written $a\cdot b\in A$ (and not $ab$, which represents an element in $A^{\otimes 2}$).

\begin{defi} 
The shuffle bialgebra $\Sh(A)=\bigoplus\limits_{n\in\mathbb{N}}\Sh_n(A) $ is the graded connected (i.e. $\Sh_0(A)=\mathbf{Q}$) commutative Hopf algebra such that
\begin{itemize}
 \item As vector spaces $\Sh_n(A)=T_n(A)$.
\item Its product $\shuffle$ is defined recursively as the sum of the two half-shuffle products $\prec ,\ \succ $
$$a_1...a_n \prec  b_1...b_p := a_1\times (a_2...a_n \shuffle  b_1...b_p)$$
$$a_1...a_n \succ  b_1...b_p := b_1...b_p \prec  a_1...a_n=b_1\times (a_1...a_n \shuffle  b_2...b_p),$$
and $\shuffle = \prec  + \succ $.
\item Its coalgebra structure is defined by the deconcatenation coproduct
$$\Delta (a_1...a_n):=\sum\limits_{0\leq k\leq n} a_1...a_k\otimes a_{k+1}...a_n.$$
\end{itemize}
\end{defi}

Recall that the notions of connected commutative Hopf algebra and connected commutative 
bialgebra are equivalent since a graded connected commutative bialgebra always has an antipode \cite{mm}. 

The above construction is functorial: $\Sh: A\rightarrow \Sh(A)$ is a functor from the category
of commutative algebras without a unit to the category of graded connected commutative Hopf algebras. 
Since a vector space can be viewed as a commutative algebra with the null product, 
our definition of $\Sh(A)$ encompasses the contruction of the shuffle algebra over a vector space. 
Although it may seem artificial at the moment (since we did not make use of the internal product $\cdot$ to define $\Sh(A)$), 
viewing $\Sh$ as a functor from commutative algebras to graded connected commutative Hopf algebras is better suited to the forthcoming developments.

Equivalently, for all $a_1,\ldots,a_{k+l} \in A$,
\begin{eqnarray*}
a_1\ldots a_k \prec a_{k+1}\ldots a_{k+l}&=& \sum_{\alpha\in \Des_{\subseteq \{k\}},\:\alpha^{-1}(1)=1}
a_{\alpha^{-1}(1)}\ldots a_{\alpha^{-1}(k+l)},\\
a_1\ldots a_k \succ a_{k+1}\ldots a_{k+l}&=& \sum_{\alpha\in \Des_{\subseteq \{k\}},\:\alpha^{-1}(1)=k+1}
a_{\alpha^{-1}(1)}\ldots a_{\alpha^{-1}(k+l)},
\end{eqnarray*}
where the $\alpha$ are permutations  of $[k+l]=\{1,...,k+l\}$.

The notation $\alpha\in \Des_{\subseteq \{k\}}$ means that $\alpha$ has at most one descent in position $k$. 
Recall that a permutation $\sigma$ of $[n]$ is said to have a descent in position $i<n$ if $\sigma(i)>\sigma(i+1)$. 
The descent set of $\sigma$, $desc(\alpha)$ is the set of all descents of $\alpha$,
$$
desc(\sigma):=\{i<n ,\sigma(i)>\sigma(i+1)\}.
$$
For $I\subset [n]$, we write $\Des_{I}:=\{\sigma , \desc(\sigma)=I\}$ and $\Des_{\subseteq I}:=\{\sigma , \desc(\sigma)\subseteq I\}.$

\medskip
A Theorem due to Sch\"utzenberger \cite{schutz} characterizes abstractly the shuffle algebras:
\begin{prop}\label{freeness}
 As a commutative algebra, $\Sh(A)$ is the \it free algebra \rm over the vector space $A$ for the relation
 \begin{equation}\label{shuffle}
(a\prec b)\prec c=a\prec (b\prec c+c\prec b).
\end{equation}
\end{prop}

An algebra equipped with a product $\prec$ satisfying this relation is sometimes refered to as a chronological algebra (although the term has also other meanings) or as a Zinbiel algebra (because it is dual to Cuvier's notion of Leibniz algebra), we refer to \cite{FoissyPatras} for historical details.
 
A last ingredient of the theory of tensors and shuffle algebras over commutative algebras will be useful: namely, the nonlinear Schur-Weyl duality established in \cite{npt}.
Let us recall that $\WQSym$ stands for the Hopf algebra of word quasi-symmetric functions. This algebra can be given various equivalent realizations (that is, its elements can be encoded by means of surjections, packed words, set compositions or
faces of permutohedra)
 and carries various algebraic structures on which we will come back later. We refer, {\it e.g.}, to \cite{npt} for details. It will appear later in this article that, whereas the realization of $\WQSym$ in terms of surjections is the one suited for studying deformations of shuffle algebras, that in terms of words is the one suited to the analysis of Hausdorff series.

For the time being, we simply recall that $\WQSym$ can be realized as the linear span of all surjections $f$ from $[n]$ to $[p]$, where $n$ runs over the integers and $1\leq p\leq n$ and postpone the definition of the product and the coproduct. We say that such a map $f$ is of degree $n$, relative degree $n-p$,  and bidegree $(n,p)$. The linear span of degree $n$ (resp. bidegree $(n,p)$) elements is written $\WQSym_n$ (resp. $\WQSym_{n,p}$). We write $\widehat{\WQSym}$ for $\prod\limits_{n,p}\WQSym_{n,p}$. 

We consider now natural
endomorphisms of the functor $T$,
viewed as a functor from nonunital commutative algebras to vector
spaces. Concretely, we look for families of linear maps $\mu_A$ from $T_n(A)$
to $ T_m(A)$ (where $A$ runs over nonunital commutative algebras and $m$
and $n$ run over the nonzero integers) such that, for any map $f$ of nonunital commutative algebras from $A$
to $B$,
\begin{equation}
\label{natu}
T_m(f)\circ \mu_A=\mu_B\circ T_n(f).
\end{equation}
Let us say that such a family $\mu_A$ satisfies the nonlinear Schur-Weyl duality
(with parameters $n,m$).
We have \cite{npt}:

\begin{prop}\label{schurweyl}
The vector space of linear maps that
satisfy the nonlinear Schur-Weyl duality with parameters $n,m$ is canonically
isomorphic to $\WQSym_{n,m}$, the linear span of surjections from $[n]$ to $[m]$.
Equivalently,
the vector space of natural endomorphisms of the functor $T$  is canonically isomorphic to
$\widehat\WQSym$.
\end{prop}

\section{Natural coalgebra endomorphisms}
The nonlinear Schur-Weyl duality shows that $\widehat\WQSym$ is the natural object for investigating 
the linear structure of the tensor and shuffle algebras over a commutative algebra.
In this section and the following ones, we study shuffle algebras from a refined point of view. 
Namely, we aim at characterizing, inside $\widehat\WQSym$, the linear endomorphisms that preserve some extra structure.
Particularly important from this point of view are the coalgebra endomorphisms, whose classification is the object of the present section.

>From now on, coalgebras $C$ are coaugmented, counital and conilpotent. 
That is, writing $\Delta$ the coproduct and $\id_C$ the identity map of $C$: 
the coalgebra is equipped with a map $\eta_C$  (the counit) from $C$ to the ground field $\Q$
which satisfies 
$$(\eta_C\otimes \id_C)\circ \Delta=\id_C=(\id_C\otimes\eta_C)\circ\Delta.$$
The ground field embeds (as a coalgebra) into $C$ so that, canonically,
$C=C^+\oplus \Q$ with $C^+:=Ker\ \eta_C$
(coaugmentation)  and, finally, writing $\overline\Delta$ the reduced coproduct 
($\overline\Delta (x) :=\Delta (x)-x\otimes 1-1\otimes x$), 
for all $x\in C^+$, there exists an integer $n$ (depending on $x$) such that 
$\overline\Delta^n(x)=0$ (conilpotency).
Coalgebra maps $\phi: C\rightarrow C'$ are required to preserve the coaugmentation and the counit, so that $\phi(C^+)\subset C'^+$. 

Recall that $\Sh(A)$ is equipped with a coalgebra structure 
$$\Delta (a_1...a_n):=\sum\limits_{i=0}^na_1...a_i\otimes a_{i+1}...a_n\in \Sh(A)\otimes \Sh(A).$$
Its counit is the projection onto $\Sh_0(A)=\Q$ and we write 
$$\Sh(A)^+=\bigoplus\limits_{n\geq 1}A^{\otimes n}.$$

The usual universal properties of the tensor algebra $T(A)$ as a free associative algebra over $A$ 
(viewed as a vector space) dualize, and we have the adjunction property
$$\Hom_{\lin}(C^+, A)\cong \Hom_{\coalg}(C,\Sh(A)),$$
where $C$ runs over coalgebras and $\Hom_{\lin}$, resp. $\Hom_{\coalg}$, stand for the set of linear maps, resp. coalgebra maps. 
In other terms, the coaugmentation coideal functor from coalgebras to vector spaces is left adjoint to the free coalgebra functor. Recall that in this statement coalgebra means coaugmented, counital and conilpotent coalgebra.
The fact that we consider only  conilpotent coalgebras is essential for the adjunction to hold, see e.g. \cite{Block} for details on the structure of cofree coalgebras in general.
In particular, we get:

\begin{cor}
There is a canonical bijection
$$\Hom_{\lin}(\Sh^+(A), A)\cong \Hom_{\coalg}(\Sh(A),\Sh(A)).$$
\end{cor}
That is, a coalgebra endomorphism $\phi$ of $\Sh(A)$ is entirely determined by the knowledge of 
$f:=\pi\circ \phi$, where we write $\pi$ for the projection from $\Sh(A)$ to $A$ orthogonally 
to $A^0=\Q$ and to the $A^{\otimes n},\ n> 1$. Conversely, any map $f\in \Hom_{\lin}(\Sh^+(A), A)$ determines a unique coalgebra endomorphism $\phi$ of $\Sh(A)$ by
$$\phi(a_1...a_n):=\sum\limits_{i_1+...+i_k=n}f(a_1...a_{i_1})\otimes \dots\otimes f(a_{i_1+....+i_{k-1}+1}...a_n).$$

A triangularity argument that we omit shows that $\phi$ is a coalgebra automorphism if and only if the restriction of $f$ to $A$ is a linear isomorphism.
For reasons that will become clear soon, we say that $f\in \Hom_{\lin}(\Sh^+(A), A)$ is tangent to identity if its restriction to $A$ (that is to a linear endomorphism of $A$) is the identity map.

Recall that by natural endomophism of the functor $\Sh$ viewed as a functor from commutative nonunital algebras to coalgebras is meant a familly $\mu_A$ (indexed by commutative nonunital algebras $A$) of coalgebra endomorphisms of the $\Sh(A)$ commuting with the morphisms (from $\Sh(A)$ to $\Sh(B)$) induced by algebra maps (from $A$ to $B$).

\begin{theo}
Let $\mathbf{Coalg}$ be the monoid of natural endomophisms of the functor $\Sh$ viewed as a functor from commutative nonunital algebras to coalgebras. 
Then, there is an isomorphism between $\mathbf{Coalg}$ and the monoid $\mathbf{Diff}$ of formal power series without  constant term, $\mathbf{Coalg}\cong X\Q[[X]]$ equipped with the substitution product (for $P(X),Q(X)\in X\Q[[X]]$, $P\circ Q(X):=P(Q(X))$.

In particular, the set $\mathbf{Coalg_1}$ of tangent-to-identity natural endomorphisms of the functor $\Sh$ is a group canonically in bijection with the group $\mathbf{Diff_1}=X+X^2\Q[[X]]$ of tangent-to-identity formal diffeomorphisms.
\end{theo} 
Let us prove first that $\mathbf{Coalg}\cong X\Q[[X]]$. Since we have a natural isomorphism 
$\Hom_{\lin}(\Sh^+(A), A)\cong \Hom_{\coalg}(\Sh(A),\Sh(A)),$ $\mathbf{Coalg}$ is canonically 
in bijection with natural transformations from $\Sh^+$ to the identity functor 
(viewed now as functors from commutative algebras to vector spaces).
By Schur-Weyl duality, we get $\mathbf{Coalg}\cong \prod\limits_{n\geq 1} \WQSym_{n,1}$. 
Identifying the unique surjection from $[n]$ to $[1]$ with the monomial $X^n$ yields the bijection $\mathbf{Coalg}\cong X\Q[[X]]$. 

We will write from now on $\phi_{P}$ for the element in $\Hom_{\coalg}(\Sh(A),\Sh(A))$ associated with a 
given formal power series $P(X)\in X\Q[[X]]$ and $f_{P}$ for the corresponding element in $\Hom_{\lin}(\Sh^+(A), A)$.

Notice that for $P(X)=\sum\limits_{i=1}^\infty p_iX^i$, the action of $f_{P}$ and $\phi_{P}$ on an arbitrary tensor $a_1...a_n\in \Sh(A)$ can be described explicitely
 \begin{equation}\label{flaf}
f_{P}(a_1...a_n)=p_n\cdot (a_1\cdot ...\cdot a_n)\in A,
\end{equation}
 \begin{equation}\label{flaphi}
\phi_{P}(a_1...a_n)=\sum\limits_{k=1}^n\sum\limits_{i_1+...+i_k=n}p_{i_1}...p_{i_k}(a_1\cdot ...\cdot a_{i_1})\otimes ...\otimes (a_{i_1+...+i_{k-1}+1}\cdot ...\cdot a_{n})
\end{equation}
This last formula describes the embedding of $\mathbf{Coalg}$ into $\widehat\WQSym$ induced by Schur-Weyl duality (the tensor product $(a_1\cdot ...\cdot a_{i_1})\otimes ...\otimes (a_{i_1+...+i_{k-1}+1}\cdot ...\cdot a_{n})$ corresponding to the nondecreasing surjection from $[n]$ to $[k]$ sending the first $i_1$ integers to $1$,..., the last $i_k$ integers to $k$).
This embedding is of course different from the one induced by the bijection with $\prod\limits_{n\geq 1} \WQSym_{n,1}$ and corresponds to the fact that elements in $\mathbf{Coalg}$ can be represented equivalently by a $f_P$ or a $\phi_{P}$: the $f_P$s are naturally encoded by elements in $\prod\limits_{n\geq 1} \WQSym_{n,1}$ (Fla (\ref{flaf})), whereas the $\phi_{P}$s are most naturally encoded by elements in $\widehat\WQSym$ (Fla (\ref{flaphi})).

Let us show now that, for arbitrary $P(X),\ Q(X)\in X\Q[[X]]$, 
\begin{equation}
\phi_{P}\circ\phi_{Q}=\phi_{P(Q)},
\end{equation}
where $P(Q)(X):=P(Q(X))$.
For an arbitrary commutative algebra $A$ and $a_1,...,a_n\in A$, we have indeed (with self-explaining notations for the coefficients of $P$ and $Q$) 
$$\pi\circ\phi_{P}\circ\phi_{Q}(a_1...a_n)=$$
$$=f_{P}(\sum\limits_{k=1}^n\sum\limits_{i_1+...+i_k=n}q_{i_1}...q_{i_k}(a_1\cdot ...\cdot a_{i_1})\otimes ...\otimes (a_{i_1+...+i_{k-1}+1}\cdot ...\cdot a_{n}))$$
$$=\sum\limits_{k=1}^n\sum\limits_{i_1+...+i_k=n}p_kq_{i_1}...q_{i_k} (a_1\cdot ...\cdot a_n)=f_{P(Q)}(a_1...a_n)=\pi\circ\phi_{P(Q)}(a_1...a_n).$$
Thus, $\phi_{P}\circ\phi_{Q}=\phi_{P(Q)}$ and the theorem follows.

\section{Formal diffeomorphisms and $\WQSym$}

We have shown that $\mathbf{Coalg}$ and $\Di$ embed naturally into $\widehat\WQSym$ (Fla (\ref{flaphi})).
We have already noticed that they are actually embedded in $\widehat\IQSym$, where $\IQSym$ stands for the linear 
span of nondecreasing surjections. Since a nondecreasing surjection $f$ from $[n]$ to $[k]$ is characterized by 
the number of elements $f_i:=|f^{-1}(i)|$ in the inverse images  of the elements of $i\in [k]$, 
nondecreasing surjections from $[n]$ to $[k]$ are in bijection with compositions of $n$ of length $k$, 
that is, ordered sequences of integers $f_1,...,f_k$ adding up to $n$. Said otherwise, nondecreasing surjections 
are naturally in bijection with a linear basis of $\Sym$ and $\QSym$, respectively the Hopf algebra of 
noncommutative symmetric functions and the dual Hopf algebra of quasi-symmetric functions, see \cite{NCSF1}, 
although the linear embeddings of $\Sym$ and $\QSym$ into $\WQSym$ induced by the bijection between the basis 
and the embedding of $\IQSym$ into $\WQSym$ 
are not standard ones. 

We shall return later on these various embeddings into $\WQSym$;
the present section studies the compatibility relations between the group structure of $\Di_1$ and the coalgebra structure existing on $\WQSym$.

Let us recall the relevant definitions. The word realization of $\WQSym$ to be introduced now will 
be useful when we will discuss later some of its combinatorial properties.
We denote by $A=\{a_1<a_2<\dots\}$ an
infinite linearly ordered alphabet and by $A^*$ the corresponding set of words.

The \emph{packed word} $u=\pack(w)$ associated with a word $w\in A^*$ is
obtained by the following process. If $b_1<b_2<\dots <b_r$ are the letters
occuring in $w$, $u$ is the image of $w$ by the homomorphism $b_i\mapsto a_i$.
For example, if $A=\N^\ast$, $pack(3\ 5\ 3\ 8\ 1)=2\ 3\ 2\ 4\ 1$.
A word $u$ is said to be \emph{packed} if $\pack(u)=u$. We denote by $\PW$ the
set of packed words.
With a word $u\in \PW$, we associate the noncommutative polynomial\footnote{
As is customary in this theory, ``polynomial'' means a formal
series of bounded degree in an infinite number of noncommuting variables, where each monomial of finite degree 
may carry a non zero scalar coefficient. We still denote by $Q\<A\>$  the corresponding algebra. 
Since these are elements of a projective limit
of polynomial rings, purists may want to call these objects {\it (noncommutative) prolynomials}.}
\begin{equation}\label{Mu}
\M_u(A) :=\sum_{\pack(w)=u}w\,.
\end{equation}
For example, restricting $A$ to the first five integers,
\begin{equation}
\begin{split}
\M_{13132}(A)=&\ \ \ 13132 + 14142 + 14143 + 24243 \\
&+ 15152 + 15153 + 25253 + 15154 + 25254 + 35354.
\end{split}
\end{equation}
Packed words $u=u_1...u_n$ are in bijection with surjections: taking for $A$ the set of integers, 
if $1,...,p$ are the letters occuring in $u$, the surjection associated with $u$ is simply the 
map from $[n]$ to $[p]$ defined by $f(i):=u_i$. 
The $\M_u$ can therefore be chosen as linear generators of $\WQSym$. Since the product of two 
$\M_u$s is a linear combination of $\M_u$s, this presentation induces an algebra structure on $\WQSym$ 
and an algebra embedding into $\Q\<A\>$. This algebra structure on $\WQSym$ is closely related to the 
Hopf algebra structure of $\QSh(A)$, see \cite{npt} for details.
Since both presentations of $\WQSym$ ($\M_u$ or surjections) are equivalent, we will not distinguish between  them, except notationally.

As for classical symmetric functions, the nature of the ordered alphabet $A$
chosen to define word quasi-symmetric functions $\M_u(A)$ is largely
irrelevant provided it has enough elements. We will therefore often omit the
$A$-dependency and write simply $\M_u$ for $\M_u(A)$, except when we want to
emphasize this dependency (and similarly for the other types of generalized
symmetric functions we will have to deal with).

The important point for us now is that
 $\WQSym$ carries naturally a Hopf algebra structure for the coproduct: 
\begin{equation}
\Delta(\M_u) :=
  \sum_{i=0}^n \M_{u_{|[1,i]}}\otimes\M_{\pack(u_{|[i+1,n]}))}.
\end{equation}
Here, $u$ is a packed word over the letters $1,...,n$ and, for
an arbitrary subset $S$ of $[n]$, $u_{|S}$ stands for the word obtained from
$u$ by erasing all the letters that do not belong to $S$.

\begin{lemma}The direct sum of the spaces of nondecreasing surjections, and of the scalars, $\IQSym$, is a subcoalgebra of $\WQSym$. 
It is a cofree coalgebra, cogenerated by the set $\Gamma\cong \N^+$ of ``elementary'' surjections $\gamma_n$ from $[n]$ to $[1]$, 
$n\geq 1$. That is, as a coalgebra, $\IQSym$ identifies with $T(\Gamma )$ (the linear span of words over the alphabet of elementary surjections) 
equipped with the deconcatenation coproduct.
\end{lemma}

In particular, $\IQSym$ is isomorphic as a coalgebra to $\QSym$, the coalgebra of quasi-symmetric functions \cite{Gessel,MR,Stan} 
and the embedding of $\IQSym$ in $\WQSym$ induces a coalgebra embedding of $\QSym$ in $\WQSym$.

Recall that coalgebra means here coaugmented counital conilpotent coalgebra, so that the free coalgebra 
over a generating set $X$ identifies with the tensor algebra over $\Q X$ equipped with the deconcatenation coproduct.
The Lemma follows then from the definition of $\Delta$. Indeed, if we write $\gamma_{i_1}...\gamma_{i_k}$ for the unique surjection 
$f$ from $[i_1+...+i_k]$ to $[k]$ such that $f(1)=...=f(i_1):=1,\ \dots ,f(i_1+...+i_{j-1}+1)=...=f(i_1+...+i_k):=k$, then,
$$\Delta(\gamma_{i_1}...\gamma_{i_k})=\sum\limits_{j=0}^k\gamma_{i_1}...\gamma_{i_j}\otimes \gamma_{i_{j+1}}...\gamma_{i_k}.$$

\begin{prop}
The embedding of $\Di_1$ into $\widehat\IQSym$ factorizes through the set of grouplike elements in $\widehat\IQSym$ (the same statement holds if we replace $\IQSym$ by the isomorphic coalgebra $\QSym$). In other terms, the group structure of $\Di_1$ is compatible with the coalgebra structure of $\WQSym$ (resp. $\QSym$).
\end{prop}

Let $P(X)=X+\sum_{i>1}p_iX^i$ and $p_1:=1$. Then, in the basis $\M_u$, 
$$\phi_{P}=\sum_{n\geq 0}\sum\limits_{k=1}^n\sum\limits_{i_1+...+i_k=n}p_{i_1}...p_{i_k}M_{1^{i_1}...k^{i_k}},$$
where $1^{i_1}...k^{i_k}$ stands for the nondecreasing packed word with $i_1$ copies of 1, ..., $i_k$ copies of $k$ 
(so that e.g. $1^32^23^2=1112233$) with the convention that the component $n=0$ of the sum contributes  $1\in\Q$ (this corresponds to $\phi_P(1)=1$).
We get
$$\Delta(\phi_{P})=\sum_{n\geq 0}\sum\limits_{k=1}^n\sum\limits_{a+b=k}\sum\limits_{i_1+...+i_k=n}p_{i_1}...p_{i_a}M_{1^{i_1}...a^{i_a}}\otimes p_{i_{a+1}}...p_{i_k}M_{1^{i_{a+1}}...(b-a)^{i_k}} $$
$$=\phi_{P}\otimes\phi_{P}.$$
 
Notice that, as a by-product, we have characterized the grouplike elements $\phi$ in $\widehat\WQSym$ which satisfy the compatibility property
\begin{equation}\label{eq:compat}
\Delta(\phi)\circ\Delta (a_1...a_n)=\Delta\circ \phi(a_1...a_n).
\end{equation}
In general, for $\phi$ in $\widehat\WQSym$ (not necessarily grouplike),  equation \eqref{eq:compat}
cannot hold if $\phi\notin \widehat\IQSym$ (this is because the deconcatenation coproduct $\Delta$ preserves 
the relative ordering of the $a_i$s). 
The converse statement is true and follows from the definition of the coproduct on $\WQSym$ (its proof is left to the reader):

\begin{lemma}Equation \eqref{eq:compat} holds for all
 $\phi\in \widehat\IQSym$,
This property characterizes $\widehat\IQSym$ as a subspace of $\widehat\WQSym$.
\end{lemma}

\begin{prop}
  The internal product (defined as the composition of surjections on $\IQSym^+$, 
the usual product on $\Q=\IQSym_0=\IQSym\cap\WQSym_0$ --the product of elements in $\IQSym^+$ and $\Q$ 
is defined to be the null product) provides $\IQSym$ with a bialgebra\footnote{
By ``bialgebra'', we simply mean in the present article a compatible product and coproduct as in  identity (\ref{equation10}), 
without requiring extra properties (very often one requires the coproduct and the product to have also compatibility 
properties with the unit of the algebra and the counit of the coalgebra).
}
structure. 
That is, the coproduct $\Delta$ and the composition of surjections $\circ$ satisfy
\begin{equation} \label{equation10} \Delta (f\circ g)=\Delta (f)\circ \Delta(g).
\end{equation}
\end{prop}

Indeed, by nonlinear Schur-Weyl duality, the identity is true if and only if the following identity holds for arbitrary $A$, and $ a_1,...,a_n$ in $A$
$$\Delta (f\circ g)\circ\Delta (a_1...a_n)=\Delta (f)\circ \Delta(g)\circ\Delta (a_1...a_n),$$
which follows from the previous Lemma and from the stability of $\IQSym$ by composition (the composition of two nondecreasing surjections is an nondecreasing surjection).

The identity map $I\in\widehat\IQSym$ (the sum of all identity maps on the finite sets $[n]$) is a unit for $\circ$ 
and behaves as a grouplike element for $\Delta$. However, $I\notin\IQSym$ and problems arise if one tries to provide $\IQSym$ with 
a classical \it graded \rm unital bialgebra structure.
A natural gradation of $\IQSym$ would be 
the relative degree of surjections (since the relative degree of the composition of two surjections $f$ and $g$ is the sum of their relative degrees). 
However, for this gradation, each graded component of $\IQSym$ is infinite dimensional (even in degree 0), 
so that many of the usual arguments regarding graded Hopf algebras do not apply directly to $\IQSym$. 

Since we aim at providing a group-theoretical picture of the theory of shuffle algebras, our final goal is to understand the 
fine structure of the image of $\Di$ in $\widehat\IQSym$.
The next section aims at clarifying these questions.

\section{Natural coderivations and the Fa\`a di Bruno algebra}

In the previous section, we have characterized the natural coalgebra endomorphisms of shuffle algebras,
or, equivalently, grouplike elements in $\widehat\IQSym$.
We have also shown that the tangent-to-identity elements form a group isomorphic to the group of tangent-to-identity formal diffeomorphisms.
We are going to study now the corresponding Lie algebra $L$ of natural coderivations of shuffle algebras.

Recall the canonical isomorphism $\IQSym\cong T(\Q\Gamma)$ (with $\Gamma\cong\N^\ast$). Let $S(\Gamma)$ stand
for the subspace of symmetric tensors in $T(\Q\Gamma)$ and $\SIQSym$ be the corresponding subspace of $\IQSym$.  
By construction, the embedding of $\Di$ into $\IQSym$ factorizes through $\SIQSym$ (see Eq. (\ref{flaphi})). 
Since symmetric tensors form a subcoalgebra of the tensor algebra for the deconcatenation coproduct 
and since the composition of two elements in $\SIQSym$ is still in $\SIQSym$,
the following Lemma is a consequence of our previous results.

\begin{lemma}
 The embedding of $\SIQSym$ into $\IQSym$ is an embedding of bialgebras (for the composition product).
\end{lemma}

Let us call tangent-to-identity an element $\mu$ in $\widehat{S(\Gamma)}$ if $\mu -I$ is a (possibly  infinite) 
linear combination of  nondecreasing strict surjections (nondecreasing surjections from $[n]$ to $[m]$ with $n>m$).

A coderivation $D$ in $\Sh(A)$ (that is, a linear endomorphism such that $\Delta\circ D=(D\otimes I+I\otimes D)\circ \Delta$, 
where $I$ stands for the identity map) is called infinitesimal if its restriction to $A\subset \Sh(A)$ 
is the null map. As usual, a natural coderivation of the shuffle algebras is a familly of coderivations 
(of the $\Sh(A)$) commuting with the morphisms induced by algebra maps 
(from $A$ to $B$, where $A$ and $B$ run over commutative nonunital algebras).

\begin{lemma}
Natural tangent-to-identity coalgebra endomorphisms of shuffle algebras identify with tangent-to-identity grouplike elements in 
$\widehat\SIQSym\cong\widehat{S(\Gamma)}$. The corresponding Lie algebra of primitive elements in $\widehat{S(\Gamma)}$ 
is the Lie algebra of natural infinitesimal coalgebra coderivations of the shuffle algebras. 
It is canonically in bijection with the Lie algebra of formal power  series $X^2\Q[[X]]$ equipped with the Lie bracket $[X^m,X^n]:=(m-n)X^{m+n-1}$.
\end{lemma}

The Lemma follows from the general property according to which, in a bialgebra, grouplike elements 
and primitive elements are in bijection through the logarithm and exponential maps, provided these maps make sense 
(that is, provided no convergence issue of the series arises). This is because, formally, for $\phi$ a grouplike element,
$$\Delta (\log(\phi))=\log (\Delta(\phi))=\log (\phi\otimes\phi)=\log(\phi)\otimes I+I\otimes\log(\phi),$$
since $\log(ab)=\log(a)+\log(b)$ when $a$ and $b$ commute and since $I$ is the unit element for the composition product.

The formal convergence of the series under consideration in the present case is ensured by the fact that a 
surjection from $n$ to $p<n$ can be written as the product of at most $n-p$ strict surjections 
(so that the coefficient of such a  surjection in the expansion of $\log(\phi)$ is necessarily finite 
and equal to its coefficient in the expansion of the truncation of the logarithmic series at order $n-p$).

The last statement of the Lemma follows from the isomorphism between tangent-to-identity coalgebra endomorphisms of 
shuffle algebras and tangent-to-identity formal diffeomorphisms. 
It can also be deduced directly from the adjunction property 
(dual to the one according to which derivations in the tensor algebra are in bijection 
with linear morphisms from $V$ to $T(V)$): writing $\Coder^+(\Sh(A))$ for the coderivations of $\Sh(A)$ vanishing on $\Q$, we have
$\Coder^+(\Sh(A))\cong \Lin(\Sh^+(A),A)$, which implies, by nonlinear Schur-Weyl duality that natural 
coalgebra coderivations of the shuffle algebras are canonically in bijection with $X\Q[[X]]$.

This bijection with $X\Q[[X]]$ can be made explicit: dualizing the formula for 
the derivation associated with
a map $f:V\longrightarrow T(V)$ 
 $$f(v_1...v_n):=\sum\limits_{i=1}^nv_1...v_{i-1}f(v_i)v_{i+1}...v_n,$$ 
we get, for $P=\sum_{i\ge1} p_iX^i$
$$D_P(a_1...a_n)=\sum_{i=1}^n\sum\limits_{j=1}^{n-i+1}p_ia_1...a_{j-1}(a_j\cdot ...\cdot a_{j+i-1})a_{j+i}...a_n.$$
In particular,  the restriction of $\phi_P$ and $D_P$ to maps from $\Sh(A)$ to $A$ agree and are both given by 
$$\phi_P(a_1...a_n)=D_P(a_1...a_n)=p_na_1\cdot ...\cdot a_n.$$

The simplest example of a coderivation is for $P(X)=X$: it is the degree operator, $Y:=D_X$,
$$Y(a_1...a_n)=\sum\limits_{i=1}^na_1...a_{i-1}I(a_i)a_{i+1}...a_n=n\ a_1...a_n.$$
Similarly, $D_{\lambda X}(a_1...a_n)=n\cdot \lambda\ (a_1...a_n).$
In general
we have, for arbitrary polynomials $P,Q$ and $\lambda\in\Q$,
$$D_{P}+\lambda D_{Q}=D_{P+\lambda Q}.$$

The description of the Lie algebra structure on $X\Q[[X]]$ induced by the isomorphism with the Lie algebra 
of natural coalgebra coderivations follows from the explicit formula for the action of coderivations 
($[D_{X^m},D_{X^n}]=(m-n)D_{X^{m+n-1}}$) but can also be deduced from the fact that composition of coalgebra 
endomorphisms is reflected in the composition of formal power series.
Recall indeed that the set of formal power series $X+X^2\Q[[X]]$ equipped with the composition product 
is the group of characters of the Fa\`a di Bruno Hopf algebra, see e.g. \cite{FGB}.

Summarizing our previous results, we get:

\begin{theo}
The Lie algebra of natural infinitesimal coalgebra coderivations of shuffle algebras is naturally isomorphic to 
(the completion of) the Lie algebra generated by the $X^n,\ n>1$, with Lie bracket $[X^m,X^n]:=(m-n)X^{m+n-1}$. 
Equivalently, it is isomorphic to (the completion of) the Lie algebra of primitive elements in 
the Hopf algebra dual to the Fa\`a di Bruno Hopf algebra.
\end{theo}

Here, ``completion'' is understood with respect to the underlying implicit grading of these 
Lie algebras and Hopf algebras (e.g. $X^n$ is naturally of degree $n-1$, see \cite{FGB} for details).

\section{Deformations of shuffles}

Let us restrict again our attention to tangent-to-identity coalgebra automorphisms of the $\Sh(A)$.
Any such $\Phi_P$ (with $P(X)-X\in X^2\Q[[X]]$) defines a natural deformation of $Sh$, that is, a new functor 
from commutative algebras to Hopf algebras
$$\Sh_P(A)=(\Sh(A),\Delta, \shuffle_P),$$
that is, $\Sh_P(A)$ identifies with $\Sh(A)$ (and with $T(A)$ equipped with the deconcatenation coproduct) 
as a coalgebra, but carries a new product defined by conjugacy
$$  x\shuffle_Py:=\phi_P(\phi_P^{-1}(x)\shuffle\phi_P^{-1}(y)).$$

For the sake of completeness, let us check explicitely the compatibity relation of the product $\shuffle_P$ and the coproduct:
since $\phi_P^{-1}$ is a coalgebra automorphism,
$$\Delta(x\shuffle_Py)=\Delta\circ \phi_P(\phi_P^{-1}(x)\shuffle\phi_P^{-1}(y))=(\phi_P\otimes\phi_P)\circ \Delta(\phi_P^{-1}(x)\shuffle\phi_P^{-1}(y))$$
$$=(\phi_P\otimes\phi_P)\circ (\Delta(\phi_P^{-1}(x))\shuffle\Delta(\phi_P^{-1}(y)))$$
$$=(\phi_P\otimes\phi_P)\circ ((\phi_P^{-1}\otimes\phi_P^{-1})\circ\Delta(x)\shuffle(\phi_P^{-1}\otimes\phi_P^{-1})\circ\Delta(y))$$
$$=\Delta(x)\shuffle_P\Delta(y).$$

\begin{defi}
The Hopf algebra $\Sh_P(A)$ is called the $P$-twisted shuffle algebra. It is isomorphic to $\Sh(A)$ as a Hopf algebra. \end{defi}

It inherits therefore all the properties of $\Sh(A)$.
The reader is referred to Reutenauer's book \cite{Reutenauer} for a systematic study of shuffle algebras.
As an algebra, $\Sh(A)$ is, for example, a free commutative algebra over a set of generators parametrized by Lyndon words.

The fundamental example of a deformation is provided by the ``$q$-exponential'' map 
$$E_q:=\sum_{n\in\N^\ast}\frac{q^{n-1}x^n}{n!}$$
which interpolates between the identity $E_0=x$ and  $E_1=e^x-1$.
The corresponding isomorphism between $\Sh(A)$ and $\Sh_{E_q}(A)$ is then given by
$$\phi_{E_q}(a_1...a_n)=\sum\limits_{\mathcal P}\frac{q^{n-k}}{P_1!...P_k!}a_{P_1}...a_{P_k},$$
where ${\mathcal P}=(P_1,...,P_k)$ runs over the nondecreasing partitions of $[n]$ 
($P_1\coprod ...\coprod P_k=[n]$ and $P_i<P_j$ if $i<j$; $a_{P_i}:=\prod\limits_{j\in P_i}a_j$ and $P_1!$ is a shortcut for $|P_1|!$.

\begin{lemma}
When $q=1$, $\Sh_{E_1}(A)$ identifies with $\QSh(A)$, the quasi-shuffle algebra over $A$, whose  product, written $\qshuffle$, is defined recursively by
$$a_1...a_n\qshuffle b_1...b_m:=$$
$$a_1(a_2...a_n\qshuffle b_1...b_m)+b_1(a_1...a_n\qshuffle b_2...b_m)+(a_1\cdot b_1)(a_2...a_n\qshuffle b_2...b_m).$$
\end{lemma}

The exponential map $\phi_{E_1}$ generalizes the Hoffman isomorphism between $\Sh(A)$ and $\QSh(A)$ 
introduced and studied in \cite{Hoffman} in the case where $A$ is a locally finite dimensional graded connected algebra.
Using a graded connected commutative algebra $A$ (instead of an arbitrary commutative algebra as in the present article), 
although a strong restriction in view of applications, has some technical advantages: it allows, for example, 
to treat directly $\QSh(A)$ as a graded connected Hopf algebra, making possible the use of structure theorems for such algebras 
(Cartier-Milnor-Moore, Leray...). 
The classical illustration of these phenomena is provided by the algebra of quasi-symmetric functions 
(the quasi-shuffle algebra over the monoid algebra of the positive integers) 
and the dual algebra of noncommutative symmetric functions: using the gradation on $\QSym$ induced by the one of the integers, 
$\N^\ast$, the exponential/logarithm transform amounts then to a mere change of basis 
(between a family of grouplike vs primitive generators) see \cite{Gessel,NCSF1,Hoffman} for details.

The proof amounts to showing that 
$$\phi_{E_1}(a_1...a_n\shuffle a_{n+1}...a_{n+m})=\phi_{E_1}(a_1...a_n)\qshuffle \phi_{E_1}(a_{n+1}...a_{n+m}).$$
Let us write $R_1...R_k$ for an arbitrary partition of $[n+m]$ such that for $1\leq i<j\leq n$ or  
$n+1\leq i<j\leq n+m$ $i\in R_p,\ j\in R_q\Rightarrow p\leq q$. 
The problem amounts to computing the coefficient of $a_{R_1}...a_{R_k}$ in the expansion of the left and right-hand sides of the equation. 
We leave to the reader the verification that only such tensors appear in these expansions. 
The coefficient is in both cases $\frac{1}{P_1!...P_k!}\frac{1}{Q_1!...Q_k!}$, 
where $P_i:=R_i\cap [n]$, $Q_i:=R_i\cap \{n+1,...,n+m\}$. This is straightforward for the right-hand side 
(the reader not familiar with quasi-shuffle products is encouraged to write down the tedious but straightforward details 
of the proof -using e.g. the recursive definition of $\qshuffle$). 
For the left-hand side, it follows from the identity ${|R_i|\choose |P_i|}\times \frac{1}{R_i!}=\frac{1}{P_i!Q_i!}$ 
and the fact that the number of words in the expansion of a shuffle product $x_1...x_l\shuffle y_1...y_k$ is $l+k\choose k$.

\bigskip

As we shall see in the sequel, even in the well-known special case of a graded algebra, 
the point of view developed in the present article is not without interest: 
being more general and conceptual than usual approaches to quasi-shuffle algebras, it allows the derivation of new insights 
on the fine structure of their operations, refining the results already obtained in \cite{npt}.

Interesting new phenomena do actually occur as soon as one considers natural linear endomorphisms of shuffle algebras.
We have already recalled from \cite{npt} that, by Schur-Weyl duality, they belong to $\WQSym$, which inherits an 
associative (convolution) product from the Hopf algebra structure of the $\Sh(A)$: 
for $f\in \WQSym_n,\ g\in \WQSym_m$ and arbitrary $a_1,...,a_{n+m}\in A$, where $A$ is an arbitrary commutative algebra,
$$f\shuffle g(a_1...a_{n+m}):=f(a_1...a_n)\shuffle g(a_{n+1}...a_{n+m}).$$

When $f$ and $g$ belong to $\FQSym$, the subset of permutations in $\WQSym$, this convolution product has a simple expression  
and defines the ``usual'' product on $\FQSym$, see e.g. \cite{MR,dht}.
Writing $\Sh_{n,m}$ for the set of $(n,m)$-shuffles (that is the elements $\sigma$ in the symmetric group $S_{n+m}$ 
of order $n+m$ such that $\sigma(1)<...<\sigma(n)$ and $\sigma(n+1)<...<\sigma(n+m)$), we get
$$f\shuffle g=\sum\limits_{\zeta\in \Sh_{n,m}}\zeta\circ (f\cdot g),$$
where $f\cdot g$ stands for the ``concatenation'' of permutations: $f\cdot g(i):=f(i)$ for $i\leq n$ and $f\cdot  g(i):=n+g(i-n)$ else.

Recall that word and free quasisymmetric functions ($\WQSym$ and $\FQSym$) carry a coproduct, defined on packed words $f$ over $k$ letters by 
$$\Delta(f):=\sum\limits_{i=0}^kf_{|\{1,...,i\}}\otimes \pack(f_{|\{i+1,...,k\}}),$$
this coproduct together with $\shuffle$ defines a graded Hopf algebra structure on both $\WQSym$ and $\FQSym$ 
(the grading is then defined by requiring a surjection from $[n]$ to $[p]$ to be of degree $n$). 
In particular, the embedding $\FQSym\subset \WQSym$ is an embedding of Hopf algebras for these structures. 
Recall however that $\shuffle$ \it is not \rm the usual product used when studying $\WQSym$, see below for details.

The following lemma is instrumental and will prove quite useful. Its proof is left to the reader.

\begin{lemma}\label{compsurj}
 For $f$ a nondecreasing surjection ($f\in \IQSym$) and $g\in\WQSym$, 
 $$\Delta(f\circ g)=\Delta (f)\circ\Delta(g).$$
\end{lemma}
Notice that the Lemma would not hold with $f$ arbitrary in $\WQSym$.

>From our previous considerations, \it any \rm formal power series in $X\Q[[X]]$ will allow us to define a new convolution product on 
$\WQSym$ associated with the P-twisted Hopf algebra structure of the $\Sh_P(A)$. 
This new product, written $\shuffle_P$ (resp. $\qshuffle$ when $P=E_1$, resp. $\shuffle_q$ when $P=E_q$) is defined by
$$\forall f\in WQSym_n ,\ g\in WQSym_m,\ $$
$$f\shuffle_Pg (a_1...a_{m+n}):=f(a_1...a_n)\shuffle_Pg(a_{n+1}...a_{n+m}),$$
(the product $f\shuffle_P g$ acts as the null map on tensors of length different from $n+m$).

The following result, although elementary, is stated as a theorem in view of its importance:

\begin{theo}
 For an arbitrary $P\in X+X^2\Q[[X]]$, the composition map $$f\longmapsto \Phi_P(f)=f_P:=\phi_P\circ f$$
 induces an isomorphism of bialgebras from $(\WQSym , \Delta, \shuffle)$ to $(\WQSym , \Delta, \shuffle_P)$. For $P=E_1$ (resp. $P=E_q$) , we will write simply $\Phi_1(f)$ (resp. $\Phi_p(f)$) for $\Phi_P(f)$.
This isomorphism is equivariant with respect to the composition product: $$\phi_P(f\circ g)=\phi_P(f)\circ g.$$
\end{theo}

The compatibility with the coproduct follows from Lemma \ref{compsurj}, from the definition of $f_P$ as the composition of $f$ with a
sum of nondecreasing surjections, and from the fact that $\Delta(\phi_P)=\phi_P\otimes\phi_P$.
On an other hand, the definition of the twisted product $\shuffle_P$ implies
$$f_P\shuffle_Pg_P=\phi_P\circ ( f\shuffle g)=(f\shuffle g)_P.$$

The following corollaries are motivated by the key role of $\FQSym$ 
and $\WQSym$ 
in the theory of noncommutative symmetric functions and their various application fields:

\begin{cor}
  Any $P\in X+X^2\Q[[X]]$ induces a Hopf algebra embedding of $(\FQSym ,\Delta, \shuffle )$ into $(\WQSym ,\Delta, \shuffle_P )$.
  This embedding is $S_n$-equivariant: for $\sigma, \beta\in S_n=\FQSym_n$, we have:
  $$(\sigma\circ\beta)_P=\sigma_P\circ \beta.$$
  \end{cor}

  \begin{cor}
For $P=E_1$, we get that $(\FQSym ,\Delta, \shuffle )$ embeds naturally into $(\WQSym ,\Delta, \qshuffle )$. As in the previous corollary, this embedding is $S_n$-equivariant.
\end{cor}

These corollaries allow to give an explicit formula for the embeddings.
Indeed, for an arbitrary $\sigma\in S_n$, we get (writing $1_n$ the identity permutation in $S_n$)
$$\sigma_P=(1_n)_P\circ \sigma.$$
Since $(1_n)_P$ is simply the component of $\phi_P$ in $\WQSym_n$, we get finally
$$\sigma_P=\sum\limits_{k=1}^n\sum\limits_{i_1+...+i_k=n}p_{i_1}...p_{i_k}1^{i_1}...k^{i_k}\circ \sigma,$$
where we write $1^{i_1}...k^{i_k}$ for the surjection sending the first $i_1$ integers to 1, ..., the integers from $i_1+...+i_{k-1}+1$ to $n$ to $k$.
For $E_q$ this formula simplifies:

\begin{lemma}\label{lemma18}
 Let $\sigma\in S_n$ and $\tau\in \WQSym_n$. We shall say that $\tau\propto \sigma$ if for all 
$1\leq i,j\leq n$, $\sigma(i)< \sigma(j)\Rightarrow \tau(i)\leq\tau(j)$. 
We also set $\tau !:=\prod\limits_{i=1}^{max(\tau )}|\tau^{-1}(\{i\})|!$ and $r(\tau)$ for the relative degree of $\tau$. 
Then, the Hopf algebra embedding $\Phi_p$ from $\FQSym$ into $\WQSym$ is given by
 $$\Phi_p(\sigma)=\sum\limits_{\tau\propto \sigma}\frac{q^{r(\tau)}\tau}{\tau !}.$$
\end{lemma}

Other consequences of the existence of such embeddings will be drawn in the sequel.

\section{Structure of twisted shuffle algebras}

The map $\phi_P$ defines an isomorphism from $\Sh(A)$ to $\Sh_P(A)$ for an arbitrary commutative algebra 
$A$ and an isomorphism between $\WQSym$ 
equipped with the shuffle product $\shuffle$ 
to $\WQSym$ 
equipped with the twisted shuffle product $\shuffle_P$.

In this section, we briefly develop the consequences of these isomorphisms and recover, among others, 
the results of \cite{npt} on projections onto the indecomposables in $\QSh(A)$.

Recall from \cite{Pat93,patJA} that $e_1:=\log^{\shuffle} (Id)$ (the logarithm of the identity map of $\Sh(A)$ 
computed using the shuffle product $\shuffle$) is a canonical section of the projection from $\Sh(A)$ to the 
indecomposables $\Sh(A)^+/(\Sh(A)^+)^2$. In particular, due to the structure theorems for graded connected Hopf algebras 
over a field of characteristic $0$ (Leray, in that particular case), $\Sh(A)$ is a free commutative algebra over the 
image of $e_1$. Equivalently, $e_1$ is the projection on the eigenspace of eigenvalue $k$ of the $k$-th Adams operation, 
that is, the $k$-th power of the identity $(Id)^{\shuffle k}$. 

These properties are clearly invariant by conjugacy and we get, since $\phi_P\circ Id\circ \phi_P^{-1}=Id$, 
the following description of $\Sh_P(A)$ as a free commutative algebra:

\begin{prop}
For an arbitrary tangent-to-identity $P$, $e_1:=\log^{\shuffle_P}(Id)$ 
(the logarithm of the identity for the $\shuffle_P$ product)
is a section of the projection from $\Sh_P(A)$ to the indecomposables $\Sh_P(A)^+/(\Sh_P(A)^+)^2$. 
In particular, $\Sh_P(A)$ is a free commutative algebra over the image of $e_1$. Equivalently, $e_1$ is 
the projection on the eigenspace associated with the eigenvalue $k$ of the $k$-th Adams operation, that is, the $k$-th power of the identity $(Id)^{\shuffle_P k}$. 
\end{prop}

The particular case of the quasi-shuffle algebra was investigated in \cite{npt}. The projection $e_1$ can then be computed explicitly. 
Recall that a surjection $f$ from $[n]$ to $[p]$ has a descent in position $i$ if and only if $f(i)>f(i+1)$. 
Let us call conjugate projection and write $\tilde f$ for the projection from $[n]$ to $[p]$ defined by: $\tilde f(i):=f(n-i)$. We have:

\begin{prop}\label{mps}
In $\WQSym$ equipped with the quasi-shuffle $\qshuffle$, 
$$e_1:=\log^{\qshuffle}(Id)
=\sum\limits_{n\geq 1}{\frac{1}{n}}\sum\limits_{I\models n} \frac{(-1)^{l(I)-1}}{{{n-1}\choose{l(I)-1}}}\sum\limits_{\Des(f)=[n]-\{i_{l(I)},...,i_{l(I)}+...+i_1\}}\tilde f,$$
where $I\models n$ means that $I=(i_1,...,i_{l(I)})$ is a composition of $n$. 
\end{prop}

To start investigating the word interpretation of $\WQSym$ and the combinatorial meaning of the formulas and results 
obtained so far, recall that, as any formula regarding $\WQSym$, this proposition can be translated into a result 
on words (and actually also into a result on Rota--Baxter algebras, due to the relationship established 
in \cite{egp} between $\WQSym$ and free Rota--Baxter algebras).

Recall that the elements of $\WQSym$ can be realized as formal sums of words over a totally ordered alphabet $X$. 
For example, the identity map $\Id\in \widehat\WQSym$ identifies under this correspondence with the formal sum of all nondecreasing words over $X$.

When the alphabet is taken to be the sequence of values of
a function from $[n-1]$ into an associative algebra,
this formal sum identifies with the $n$-th value of the function (unique, formal) solution to the recursion $F=1+S(F\cdot f)$, 
where $S$ is the summation operator $S(f)(j):=\sum_{i=1}^{j-1}f(i)$, $S(f)(1):=0$. 
That is, $$F(k)=1+\sum\limits_{i=1}^{k-1}\sum\limits_{1\leq a_1<...<a_k\leq k-1}f(a_1)...f(a_k).$$
Since the concatenation product of words induces an associative product on $\WQSym$ that identifies with $\qshuffle$ \cite{dht}, 
we get finally that $e_1$ computes in that case $\log(F)(n)$. This phenomenon was studied recently in detail in \cite{em2013}, to which we refer for details.

\section{Gradations}
In the context of shuffle algebras, the conjugacy map by $\phi_P$, for an arbitrary tangent-to-identity $P$, 
maps the degree operator $Y$ on $\Sh(A)$ to a degree operator $Y_P:=\phi_P\circ Y\circ \phi_P^{-1}$ on $\Sh_P(A)$.
That is, more explicitly:
\begin{prop}
The operator $Y_P:=\phi_P\circ Y\circ \phi_P^{-1}$ acting on $\Sh_P(A)$ is a derivation and a coderivation which 
leaves invariant the subspaces $\Sh_P^{\leq n}(A):=\bigoplus\limits_{i\leq n}A^{\otimes n}$. 
Its action is diagonalizable, with eigenvalues $i\in\N$. 
The eigenspaces for the eigenvalues $0$ and $1$ are the scalars, resp. $A$. 
In general, the eigenspace for the eigenvalue $n$ is contained in $\Sh_P^{\leq n}(A)$ and its 
intersection with $\Sh_P^{\leq n-1}(A)$ is the null vector space, more precisely:
$$\Sh_P^{\leq n}(A)=\Sh_P^{\leq n-1}(A) \oplus Ker(Y_p-nId).$$
\end{prop}

Concretely, conjugacy by $\phi_P$ defines an isomorphism between $A^{\otimes n}\subset \Sh(A)$ 
and the eigenspace associated with the eigenvalue $n$ of $Y_P$ in $\Sh_P(A)$. 
The conjugacy map can be described explicitely as follows:

\begin{prop}
For $U\in X+X^2\Q[[X]]$ (or more generally in $\Q^\ast X+X^2\Q[[X]]$) and $V$ an arbitrary formal power series without constant term,
$$\phi_U^{-1}\circ D_V\circ \phi_U=D_{\frac{V\circ U}{U'}}.$$
\end{prop}

By linearity of $D$  and (formal) continuity of the action by conjugacy, it is enough to prove the formula when $V=X^p$.
We denote by $W$ the inverse of $U$ for the composition and write $u_i$ the coefficients of $U$, and similarly for $V$ and $W$. Then,

$$\pi\circ\phi_W\circ D_{X^p}\circ \phi_U(a_1...a_n)$$
$$=f_W\circ D_{X^p}(\sum\limits_{k=1}^n\sum\limits_{i_1+...+i_k=n}u_{i_1}...u_{i_k}(a_1\cdot ...\cdot a_{i_1})...(a_{i_1+...+i_{k-1}+1}\cdot ...\cdot a_{n}))$$
$$=\sum\limits_{k=p}^n\sum\limits_{i_1+...+i_k=n}(k-p+1)w_{k-p+1}u_{i_1}...u_{i_k}(a_1\cdot ...\cdot a_n),$$
so that $\pi\circ\phi_W\circ D_{X^p}\circ \phi_U$ is the linear map associated to the formal power series:
$$\sum\limits_{k=p}^\infty (k-p+1)w_{k-p+1}U^k=\big( \sum\limits_{k=p}^\infty (k-p+1)w_{k-p+1}X^k\big) \circ U$$
$$=\big(\sum\limits_{i=1}^\infty iw_iX^{i-1+p}\big)\circ U=(X^pW')\circ U=U^p\cdot (W'\circ U)=\frac{U^p}{U'}.$$
Hence, $$\phi_U^{-1}\circ D_{X^p}\circ \phi_U=D_{\frac{U^p}{U'}},$$
and the proposition follows.

\ \\

In particular, taking $P=E_1$, we get:

\begin{prop}
The eigenspaces of the coderivation $D_{(1+X)\ln(1+X)}$ provide the quasi-shuffle algebras $\QSh(A)$ with a grading and, more precisely, provide the triple $(\QSh(A),\qshuffle,\Delta)$ with the structure of graded connected commutative Hopf algebra.
\end{prop}

Indeed, with $D=\psi\circ D_X\circ \psi^{-1}=\psi\circ Y\circ \psi^{-1}$ and $\psi:=\phi_{\exp(X)-1}$, 
$$D=\phi_{\ln(1+X)}^{-1}\circ D_X\circ\phi_{\ln(1+X)}=D_{(1+X)\ln(1+X)}.$$

Notice that, since $(1+X)\ln(1+X)=1+\sum\limits_{k+2}^\infty\frac{(-1)^k}{k(k+1)}X^k$,
$$D_{(1+X)\ln(1+X)}(a_1...a_n)=n\cdot a_1...a_n+\sum\limits_{i=2}^n\sum\limits_{j=1}^{n-i+1}\frac{(-1)^i}{i(i-1)}a_1...a_{j-1}(a_j\cdot ...\cdot a_{j+i-1})a_{j+i}...a_n.$$

\section{A new realization of $\Sym$}

In the last sections, we develop some combinatorial applications of our previous results. 
We will focus mainly on the consequences of Lemma \ref{lemma18}, that is, the existence 
of a new isomorphic embedding of $\FQSym$ 
into $\WQSym$.

The present section explains briefly how these results translate in terms of \it polynomial realizations \rm 
of these algebras. In particular, we emphasize that the previous embedding gives rise to 
\it new \rm realizations of $\Sym$ and $\FQSym$ that will appear in the forthcoming section 
to be meaningful for the combinatorial study of the Hausdorff series. 

Recall that $\Sym$, the Hopf algebra of noncommutative symmetric functions, is the free associative algebra 
generated by a sequence $S_i,\ i\in\N^\ast$ of divided powers ($\Delta(S_n):=\sum\limits_{i=0}^nS_i\otimes S_{n-i}$, where $S_0:=1\in \Q$).
In \cite{NCSF1}, it has been shown that this Hopf algebra could bring considerable
simplifications in the analysis of the so-called continuous BCH
(Baker-Campbell-Hausdorff) series, the formal series $\Omega(t)=\log X(t)$
expressing the logarithm of the solution of the (noncommutative) differential
equation $X'(t)=X(t)A(t)$ ($X(0)=1$) as iterated integrals of products
of factors $A(t_i)$.
The new polynomial realization of $\Sym$ to be introduced will lead instead to a
straightforward proof of Goldberg's formula
for the coefficients of the (usual) Hausdorff series.

This algebra $\Sym$ can be embedded as a  Hopf subalgebra into $\FQSym$ 
by sending $S_n$ to the identity element in the symmetric group of order $n$. This results into the usual 
polynomial realization of $\Sym$ introduced in \cite{NCSF1} 
(that is, its realization through an embedding into the algebra of noncommutative polynomials $\Q \<X\>$) 
sending $S_n$ to the sum of all nondecreasing words of length $n$. We refer e.g. to \cite{nt} for details.

Instead of doing so, we may now take advantage of Lemma \ref{lemma18} and define a new polynomial realization 
of $\Sym$ using the existence of an isomorphic embedding $\Phi_1$ of $\FQSym$ 
into $\WQSym$: 
\begin{equation}
\hat{S}_n = \sum_{u\ \text{nondecreasing}, \ |u|=n}\frac1{u!}\M_u,
\end{equation}
where for notational convenience we write $\hat{S_n}$ for $\Phi_1(S_n)$ (and similarly for the images by $\Phi_1$ of the other elements of $\Sym$ and $\FQSym$ in $\WQSym$).
In view of Fla (\ref{Mu}), we have equivalently
\begin{equation}
\hat{\sigma}_t := \sum_{n\ge 0}t^n \hat{S}_n =e^{tx_1}e^{tx_2}\cdots =\prod_{i\ge 1}^\rightarrow e^{tx_i}\,.
\end{equation}
Notice that $\hat{\sigma}_t$ is (as expected) a grouplike element for the
standard coproduct of noncommutative polynomials for which letters $x_i$ are primitive. 

An interesting feature of this realization is that $\Phi:=\log \hat{\sigma}_1$
is now the Hausdorff series
\begin{equation}
\Phi = \log(e^{x_1}e^{x_2}\cdots) = H(x_1,x_2,x_3,\dots)\,.
\end{equation}
Moreover, two nondecreasing
words $v$ and $w$ such that $\pack(v)=\pack(w)=u$ have the same
coefficient in $\hat{\sigma}_1$, that is,
\begin{equation}
\frac{1}{u!},\quad\text{where $u! := \prod_i m_i(u)!$}
\end{equation}
and $m_i(u)$ is the number of occurences of $i$ in $u$.

The Hausdorff series can now be expanded in the basis $\M_u$ of $\WQSym$ as
\begin{equation}
\Phi = \sum_{u} c_u\M_u
\end{equation}
and one may ask whether the previous formalism can shed any light on the coeffients
$c_u$. 
There is actually a formula for $c_u$, due to Goldberg \cite{Go}, and reproduced
in Reutenauer's book \cite[Th. 3.11 p. 63]{Reutenauer}. 
This formula, which was obtained
as a combinatorial {\it tour de force}, will be shown in the sequel to be a direct consequence of our previous results.

\section{Goldberg's formula revisited}

Let us fix first some notations.
For $I=(i_1,...,i_r)$, we set $S^I:=S_{i_1}...S_{i_r}$; we also write $\ell(I)=r$ and $|I|=i_1+...+i_r$ (so that $I\vDash |I|$).
By definition,
\begin{equation}
\Phi=\log(1+\hat{S}_1+\hat{S}_2+\cdots) = \sum_{r\ge 1}\frac{(-1)^{r-1}}{r} \sum_{\ell(I)=r}\hat{S}^I\\ 
= \int_{-1}^0\left(
               \sum_I t^{\ell(I)} \hat{S}^I \right)\frac{dt}{t}
\end {equation}
so that the coefficient $c_u$ of $\M_u$ in the Hausdorff series
is,
denoting by $\NN_u$ the dual basis of $\M_u$,
\begin{equation}
c_u = \int_{-1}^0\left\< \NN_u,  \sum_I t^{\ell(I)}\hat{S}^I  \right\> \frac{dt}{t}\,.
\end{equation}

For $u$ a word of length $n$, we have therefore to evaluate $\left\< \NN_u,\hat{A}_n(t)\right\>$ with $A_n(t):=\sum\limits_{|I|=n}t^{l(I)}S^I$.
This last sum  is related to a well-known series. 
The {\em noncommutative Eulerian polynomials} are defined by \cite[Section 5.4]{NCSF1}
\begin{equation}
{\mathcal A}_n(t) 
=\ \sum_{k=1}^n (\sum_{{\scriptstyle |I|=n}\atop{\scriptstyle \ell(I)=k}} {R}_I) t^k=
\ \sum_{k=1}^n \ {\bf A}(n,k)\, t^k.
\end{equation}
where $R_I$ is the ribbon basis  (the basis of $\Sym$ obtained from the $S^I$ basis by M\"obius inversion in the boolean lattice) \cite[Section 3.2]{NCSF1}.
The generating series of the ${\mathcal A}_n(t)$ is given by
\begin{equation}
{\mathcal A}(t) := \ \sum_{n\ge 0} \, {\mathcal A}_n(t)
=
(1-t) \, \left( 1 - t\, \sigma_{1-t} \right)^{-1} \,,
\end{equation}
where $\sigma_{1-t}=\sum (1-t)^nS_n$.
Let ${\mathcal A}_n^*(t) = (1-t)^{-n}\, {\mathcal A}_n(t)$.
Then,
\begin{equation}
{\mathcal A}^*(t)
:=
\ \sum_{n\ge 0} \, {\mathcal A}_n^*(t)
=
\sum_{I} \
\left( \displaystyle {t \over 1-t} \right)^{\ell(I)} \, S^I \ .
\end{equation}
and
\begin{equation} \label{GEN*}
\sum_{I\vDash n} t^{\ell(I)}S^I =A_n(t)= {\mathcal A}_n^*\left(\frac{t}{1+t}\right)
=(1+t)^n{\mathcal A}_n\left(\frac{t}{1+t}\right).
\end{equation}

To evaluate, for a packed word $u$ of length $n$, the pairing
$\<\NN_u,\hat{A}_n(t)\>$,
let us start with the observation that, if $u=1^n$, then, writing $\F_\sigma$ for the dual basis to the $\sigma\in S_n=\FQSym_n$, 
\begin{equation}
\Phi_1^\dagger(\NN_u)=\frac1{n!}\sum_{\sigma\in\SG_n}\F_\sigma,
\end{equation}
where $\Phi_1^\dagger$ is the adjoint map,
so that in this case,
\begin{equation}
\<\NN_u,\hat{A}_n(t)\> = \<\Phi_1^\dagger(\NN_u),A_n(t)\>
=\frac1{n!}\sum_{\sigma\in\SG_n}t^{d(\sigma)+1}(1+t)^{r(\sigma)}
= \frac1{n!}t E_n(t,t+1)
\end{equation}
where $d(\sigma)$ is the number of descents of $\sigma$,
$r(\sigma)=n-d(\sigma)$ the number of rises, and $E_n$ is the homogeneous
Eulerian polynomial normalized as in \cite{Reutenauer}
\begin{equation}
E_n(x,y) = \sum_{\sigma\in\SG_n}x^{d(\sigma)}y^{r(\sigma)}\,.
\end{equation} 

Now, recall that the coproduct of $\NN_u$ dual to the product of the $\M_u$ is \cite{NTtri}
\begin{equation}
\Delta\NN_u = \sum_{u=u_1u_2}\NN_{\pack(u_1)}\otimes \NN_{\pack(u_2)}
\end{equation}
(deconcatenation). We can omit the packing operation in this formula
if we make the convention that $\NN_w=\NN_u$ if $u=\pack(w)$. Then,
since $\Phi_1$, and hence also $\Phi_1^\dagger$ are morphisms of
Hopf algebras,
 for a composition
$L=(l_1,\ldots,l_p)$,
\begin{equation}
\<\Phi_1^\dagger(\NN_u),S^L\>=\<\Delta^{[k]}(\NN_u),S_{l_1}\otimes ...\otimes S_{l_p}\>=
\prod_{k=1}^p\<\Phi_1^\dagger(\NN_{u_k}),S_{l_k}\>
\end{equation}
where $\Delta^{[k]}$ is the $k$-th iterated coproduct and $u=u_1u_2\cdots u_p$ with $|u_k|=l_k$ for all $k$. Moreover, this is nonzero
only if all the $u_k$ are nondecreasing, in which case the result is
$1/(u_1!\cdots u_p!)$.

Thus, if
\begin{equation}
u = w_1\cdots w_m
\end{equation}
is the factorization of $u$ into maximal nondecreasing words,
with $|w_k|=n_k$,
we have
\begin{equation}
\<\Phi_1^\dagger(\NN_u),A_n(t)\>= \prod_{k=1}^m\<\Phi_1^\dagger(\NN_{w_k}),A_{n_k}(t)\>
\end{equation}
since 
\begin{equation}
 \prod_{k=1}^m A_{n_k}(t) = \sum_{I\in C_u} t^{\ell(I)}S^I
\end{equation}
where $C_u$ is the set of compositions which are a refinement of $(n_1,...,n_m)$ and are the ones such that 
$\<\Phi_1^\dagger(\NN_u),S^I\>\not=0$.

Next, if
$v=1^{l_1}2^{l_2}\cdots p^{l_p}$,
\begin{equation}
\<\Phi_1^\dagger(\NN_v),S^L\>=\prod_{k=1}^p\<\Phi_1^\dagger(\NN_{k^{l_k}}),S^{l_k}\>
\end{equation}
(both sides are equal to $1/(l_1!\cdots l_p!)$),
so that
\begin{equation}\label{eqvarphi}
\<\Phi_1^\dagger(\NN_v),A_{|L|}(t)\>
=\left(1+\frac1t\right)^{r(v)}
\prod_{k=1}^p\<\Phi_1^\dagger(\NN_{k^{l_k}}),A_{l_k}(t)\>\,.
\end{equation}
where $r(v)$ is the number of different letters (or of strict rises) of $v$.
Indeed, $A_{|L|}(t)=\sum\limits_{|I|=|L|}t^{l(I)}S^I$ and, since the $S_l$ are grouplike,
$$\<\Phi_1^\dagger(\NN_u),S^I\>=\prod\limits_{k=1}^p\<\Phi_1^\dagger(\NN_k^{l_k}),S^{I|k}\>,$$
where $I|k$ is the partition of $\{l_1+...+l_{k-1}+1,...,l_1+....+l_k\}$ induced by the partition $I$ of $[|L|]$. Finally, writing $I\cup L$ for the partition refining $I$ and $L$ 
(obtained, {\it e.g.}, by gluing the $I|k$), using
\begin{equation}
\<\Phi_1^\dagger(\NN_u),t^{l(I)}S^I\>
=t^{l(I)-l(I\cap L)}
\<\Phi_1^\dagger(\NN_u),t^{l(I\cap L)}S^{I\cap L}\>
\end{equation}
and noting that 
$|\{I\models |L|,\ I\cap L=K, \ l(I)-l(K)=k<r(u)\}|={r(u)\choose k}$, we get (\ref{eqvarphi}).

We can  now see that if we decompose $u$ into maximal blocks of 
identical letters,
\begin{equation}
u = i_1^{j_1}i_2^{j_2}\cdots i_s^{j_s}
\end{equation}
we have finally
$$
\<\Phi_1^\dagger(\NN_u),A_n(t)\>=
\left(1+\frac1t \right)^{r(u)}\prod_{k=1}^s\<\Phi_1^\dagger (\NN_{i_k^{j_k}}),A_{j_k}(t)\>$$
$$=t^{d(u)+1}(1+t)^{r(u)}\prod_{k=1}^s
\frac{E_{j_k}(t,1+t)}{j_k!}$$
which implies Goldberg's formula:
\begin{theo}
The coefficient $c_u$ of $\M_u$ in the Hausdorff series $\Phi$ is given by:
\begin{equation}
c_u=\int_{-1}^0
t^{d(u)+1}(1+t)^{r(u)}\prod_{k=1}^s
\frac{E_{j_k}(t,1+t)}{j_k!}\frac{dt}{t}
\end{equation}
 \end{theo}

More generally, for an arbitrary
moment generating function
\begin{equation}
f(z) = \sum_{n\ge 1}f_n z^n.
\end{equation}
with
\begin{equation}
f_n =\int_\R t^n d\mu(t)
\end{equation}
 the coefficient of $\M_u$ in $f(\hat{\sigma}_1)$ is
\begin{equation}
\int_\R t^{d(u)+1}(1+t)^{r(u)}\prod_{k=1}^s
\frac{E_{j_k}(t,1+t)}{j_k!} d\mu(t).
\end{equation}

%%%%%%%%%%%%%%%%%%%%%%%%%%%%%%%%%%

\end{document}